\documentclass[12pt,reqno]{amsart}
\usepackage{amssymb}
\usepackage{cases}
\usepackage{bbm}
\usepackage{mathrsfs}
\usepackage{graphicx}
\usepackage{amsfonts}
\usepackage{enumerate}

\usepackage{amssymb, amsmath, amsfonts, latexsym}
\usepackage{enumerate}

\newtheorem{thm}{Theorem}[section]
\newtheorem{corollary}[thm]{Corollary}
\newtheorem{lemma}[thm]{Lemma}
\newtheorem{proposition}[thm]{Proposition}

\newtheorem{remarks}[thm]{Remark}

\newtheorem{hyp}[thm]{Hypothesis}

\theoremstyle{definition}
\newtheorem{defn}{Definition}[section]
 \theoremstyle{remark}

\numberwithin{equation}{section}

\newcommand{\HS}{\mathrm{HS}}

\newcommand{\ee}{\mathbb{E}}

\newcommand{\nn}{\mathbb{N}}
\newcommand{\rr}{\mathbb{R}}
\newcommand{\pp}{\mathbb{P}}
\newcommand{\qq}{\mathbb{Q}}

\newcommand{\e}{\varepsilon}

\def\AA{\mathcal A}
\def\BB{\mathcal B}

\def\FF{\mathcal F}
\def\EE{\mathcal E}

\def\AA{\mathcal A}
\def\BB{\mathcal B}

\def\FF{\mathcal F}
\def\EE{\mathcal E}
\def\GG{\mathcal G}
\def\HH{\mathcal H}

\def\<{\langle}
\def\>{\rangle}

\def\beq{\begin{equation}}
\def\nneq{\end{equation}}

\def\bdef{\begin{defn}}
\def\ndef{\end{defn}}

\def\bthm{\begin{thm}}
\def\nthm{\end{thm}}

\def\bprop{\begin{proposition}}
\def\nprop{\end{proposition}}

\def\brmk{\begin{remarks}}
\def\nrmk{\end{remarks}}

\def\bexa{\begin{exa}}
\def\nexa{\end{exa}}

\def\blem{\begin{lemma}}
\def\nlem{\end{lemma}}

\def\bcor{\begin{corollary}}
\def\ncor{\end{corollary}}

\def\bhyp{\begin{hyp}}

\def\nhyp{\end{hyp}}

\def\<{\langle}
\def\>{\rangle}

\date{}

\def\bexe{\begin{exe}}
\def\nexe{\end{exe}}

\def\bprf{\begin{proof}}
\def\nprf{\end{proof}}

\def\bdes{\begin{description}}
\def\ndes{\end{description}}

\title[ MDP for the Langevin equation]{Moderate deviations for the Langevin equation with strong damping}

\author{Lingyan Cheng }
\address{Lingyan Cheng \\
Center of Applied Mathematics, Tianjin University, Tianjin 300072, PR China}
\email{lycheng@tju.edu.cn}

\author{Ruinan Li}
\address{Ruinan Li \\
School of Statistics and Information, Shanghai University of International Business and Economics, Shanghai 201620,  PR China}
\email{ruinanli@amss.ac.cn}

\author{Wei Liu}
\address{Wei Liu\\School of Mathematics and Statistics, Wuhan University, Wuhan, Hubei 430072, PR China;
Computational Science Hubei Key Laboratory, Wuhan University, Wuhan, Hubei 430072, PR China}
\email{wliu.math@whu.edu.cn}

\date{}

\begin{document}

\maketitle

 \noindent {\bf Abstract:}
In this paper, we establish a moderate deviation principle for the Langevin dynamics  with  strong damping.     The weak convergence approach  plays an important role in the proof.
 \vskip0.3cm

 \vskip0.3cm

\vskip0.3cm

\section{Introduction}

\noindent
For every $\e >0$, consider the following Langevin equation with strong damping
\beq\label{Langevin eq}
\begin{cases}
\ddot{q}^\e(t)=b(q^\e(t))-\frac{\alpha(q^\e(t))}{\e} \dot{q}^\e(t)+\sigma(q^\e(t))\dot{B}(t), \\
q^\e(0)=q \in \rr^d,\quad \dot{q}^\e(0)=p \in \rr^d.
\end{cases}
\nneq
Here $B(t)$ is a $d$-dimensional standard Wiener process, defined on some complete stochastic basis $(\Omega,\FF,\{\FF_t\}_{t\ge 0 },\pp)$. The coefficients $b, \alpha$ and $\sigma$ satisfy some regularity conditions (see Section 2 for details) such that for any fixed $\e>0, T>0$ and $ k \ge 1$,  Eq.\eqref{Langevin eq} admits a unique solution $q^\e$ in $L^k(\Omega;C([0,T];\rr^d))$. Let $q_\e (t):= q^\e (t/\e)$, $ t \ge 0$, then Eq.\eqref{Langevin eq} becomes
\beq\label{Langevin eq time}
\begin{cases}
\e^2 \ddot{q}_\e (t)=b(q_\e(t))-\alpha(q_\e(t)) \dot{q}_\e (t) +\sqrt{\e}\sigma (q_\e(t)) \dot{w}(t), \\
q_\e(0)=q \in \rr^d,\quad \dot{q}_\e(0)=\frac {p}{\e} \in \rr^d,
\end{cases}
\nneq
 where $w(t):=\sqrt{\e} B(t/\e)$, $ t \ge 0$ is  also a $\rr^d$-valued Wiener process.

 \vskip0.3cm

In \cite{CF},  Cerrai and Freidlin established a large deviation principle (LDP for short) for Eq.\eqref{Langevin eq time} as $\e\to 0+$.
More precisely, for any  $T>0$, they proved that  the  family $\{q_\e\}_{\e>0}$  satisfies the  LDP in the space $C([0,T]; \rr^d)$, with the same rate function $I$ and the same speed
function $\e^{-1}$ that describe the LDP of the first order equation
\beq\label{eq g}
\dot{g}_\e (t)= \frac{b(g_\e(t))}{\alpha(g_\e(t))} +\sqrt{\e} \frac{\sigma(g_\e (t))}{ \alpha(g_\e(t))} \dot{w}(t),\ \ \ \
g_\e (0) =q \in \rr^d.
\nneq
Explicitly, this means that
\begin{itemize}
\item[(1)] for any constant $c>0$, the level set $\{f; I(f)\le c\}$ is compact in $C([0,T];\rr^d)$;
\item[(2)]for any closed subset $F\subset C([0,T];\rr^d)$,
$$
\limsup_{\e\rightarrow0+}\e\log\pp(q_{\e}\in F)\le -\inf_{f\in F}I(f);
$$
\item[(3)] for any open subset $G\subset C([0,T];\rr^d)$,
$$
\liminf_{\e\rightarrow0+}\e\log\pp(q_{\e}\in G)\ge -\inf_{f\in G}I(f).
$$
\end{itemize}

  The  dynamics system \eqref{eq g} can be regarded as  the random perturbation of the following deterministic differential equation
 \beq\label{deterministic eq}
\dot{q}_0(t) = \frac{b(q_0(t))}{\alpha(q_0(t))},\ \ \
q_0(0)=q \in \rr^d.
\nneq
Roughly speaking, the LDP result in \cite{CF}  shows that  the  asymptotic probability of
$\pp(\|q_{\e} -q_0\| \ge \delta)$ converges exponentially to $0$ as $\e \to 0$ for any $\delta>0$, where $\|\cdot\|$ is the sup-norm on $C([0,T];\rr^d)$.

\vskip0.3cm

Similarly to the large deviations, the moderate deviations arise in  the theory of statistical  inference quite naturally. The  moderate deviation principle (MDP for short)   can provide us with the rate
  of convergence and a useful method for constructing asymptotic confidence intervals (see, e.g., recent works \cite{GZ, HS, Kal, MiaoShen} and references therein). Usually, the quadratic form of the rate function corresponding to the MDP
allows for the explicit minimization, and particularly it allows one
to obtain an asymptotic evaluation for the exit time (see \cite{Klebaner1999}). Recently, the study of the MDP estimates for  stochastic (partial) differential equation has been carried out as well, see   e.g. \cite{BDG, GL, WZZ, WZ} and so on.

\vskip0.3cm
 In this paper, we shall investigate the MDP problem for the family $\{q_\e\}_{\e>0}$ on $ C([0,T];\rr^d)$.  That is, the asymptotic behavior of the trajectory
\beq\label{eq X e}
X_\e (t) = \frac{1}{\sqrt{\e} h(\e)} \left(q_\e(t)-q_{0}(t)\right),\quad  t \in [0,T].
\nneq
 Here  the deviation scale satisfies
\beq\label{h}
h(\e)\to +\infty\ \text{ and }\ \sqrt{\e} h(\e) \to 0, \quad \text{as}\  \e \to 0.
\nneq

Due to the complexity of $q_\e$, we mainly use the weak convergence approach to deal with this problem. Comparing with the approximating method used in Gao and Wang \cite{GW}, our method is simpler since we only need the moment estimation rather than the exponential moment estimation of the solution.

\vskip0.3cm

 The organization of this paper is as follows. In Sect. 2, we present the framework of the Langevin equation, and then state  our main results. Sect. 3 is devoted to proving the MDP.

\section{Framework and main results}

\vskip0.3cm

Let  $|\cdot|$ be the Euclidean norm of a vector in $\rr^d$, $\langle \cdot, \cdot\rangle$ the inner production in $\rr^d$, and  $\| \cdot \|_{\HS}$ the Hilbert-Schmidt norm in $\rr^{d\times d}$ (the space of $d\times d $ matrices). For a function $b: \rr^d\rightarrow\rr^d$, $Db=\left(\frac{\partial}{\partial x_j} b^i\right)_{1\le i,j \le d}$  is the Jacobian matrix of $b$. Recall that $\|\cdot\|$ is the sup-norm on $C([0,T];\rr^d)$.
Throughout this paper, $T>0$ is some fixed constant, $C(\cdot)$ is a positive constant depending on the parameters in the  bracket and independent of $\e$. The value of $C(\cdot)$ may be different from line to line.

 Assume that the coefficients $b,\alpha$ and $\sigma$ in  \eqref{Langevin eq time} satisfy the following hypothesis.
\bhyp\label{hyp}
\begin{itemize}
 \item[(a)] The mappings $b: \rr^d \to \rr^d $ and $\sigma: \rr^d \to \rr^{d\times d}$ are continuously differentiable,  and there exists some constant $K>0$ such that for all $x, y\in \rr^d$,
    \beq\label{bLip}
    |b(x)-b(y)| \le K|x-y|,
    \nneq
    and
    $$ \|\sigma(x)-\sigma(y)\|_{\HS} \le K|x-y|,\  \|\sigma(x)\|_{\HS} \le K.
    $$
    Moreover, the matrix $\sigma(q)$ is invertible for any $q \in \rr^d$, and $\sigma^{-1}: \rr^d \to \rr^{d\times d}$ is bounded.
\item[(b)] The mapping $\alpha: \rr^d \to \rr$ belongs to $C_b^1(\rr^d)$ %(the family of bounded and continuously differentiable functions)
     and there exist  some constants $0<\alpha_0\le\alpha_1$ and $K>0$ such that
    $$
    \alpha_0=\inf_{x \in \rr^d} \alpha (x), \ \alpha_1=\sup_{x \in \rr^d}\alpha(x) \text{ and }  \sup_{x \in \rr^d}|\nabla \alpha (x)|\le K.
    $$
    \end{itemize}

\nhyp

Notice that:
\begin{itemize}
\item[(1)]
$\|Db\|_{\HS}\le K$ since $b$ is continuously differentiable and satisfies \eqref{bLip};
\item[(2)]
$\sigma/\alpha$ is Lipschitz continuous and bounded due to the Lipschitz-continuity and the boundness of the functions $\sigma$ and $1/\alpha$.
\end{itemize}

\vskip0.3cm
Under Hypothesis \ref{hyp}, according to \cite[Theorem 2.2]{GW}, we know that the family $\left\{(g_\e-q_{0})/[\sqrt{\e}h(\e)]\right\}_{\e>0}$ satisfies the LDP on $C([0,T]; \rr^d)$ with speed $h^2 (\e)$ and a good rate function $I$ given by
\beq\label{rate function 1}
I(\psi) = \frac {1}{2} \inf_{h \in \HH ; \psi=\Gamma_0(h)} \| h\|_{\HH}^2,
\nneq
where
\beq\label{HH}
\HH: =\left\{ h \in C([0,T];\rr^d);\ h(t)=\int_0^t \dot{h}(s)ds,\ \|h\|_\HH^2:= \int_0^T \left|\dot{h}(t)\right|^2 dt <\infty \right\}
\nneq
 and
\beq\label{Gamma}
\Gamma_0(h(t))= \int_0^t D\left( \frac{b(q_0(s))}{\alpha(q_0(s))}\right) \Gamma_0(h(s))ds +\int_0^t \frac{\sigma (q_0(s))}{\alpha(q_0(s))} \dot{h}(s) ds,
\nneq
with the convention $\inf\emptyset=\infty$. This special kind of LDP is just the MDP for the family $\{g_\e\}_{\e>0}$ (see \cite{DZ}).

\vskip0.3cm
The main goal of this paper is to prove that the family $\{q_\e\}_{\e>0}$ satisfies  the same MDP as  the family $\{g_\e\}_{\e>0}$ on $ C([0,T];\rr^d)$.
\bthm\label{MDP}
Under Hypothesis \ref{hyp}, the family $\{(q_\e-q_{0})/[\sqrt{\e}h(\e)] \}_{\e >0}$ obeys an LDP on $C([0,T]; \rr^d)$ with the speed function $h^2(\e)$ and the rate function $I$ given by \eqref{rate function 1}.
\nthm

\section{ Proof of MDP }

\subsection{Weak convergence approach in LDP}

In this subsection, we will give the general criteria for the LDP given in \cite{BDM}.
\vskip0.3cm
Let $(\Omega,\mathcal{F},\mathbb{P})$ be a probability space with an increasing family $\{\FF_t\}_{0\le t\le T}$ of the sub-$\sigma$-fields of $\FF$ satisfying the usual conditions.
Let $\mathcal{E}$ be a Polish space with the Borel $\sigma$-field $\mathcal{B}(\mathcal{E})$.
The Cameron-Martin space associated with the Wiener process $\{w(t)\}_{0\le t\le T}$ (defined on the filtered probability space given above) is given by \eqref{HH}.
See \cite{DZ}. The space $\HH$ is a Hilbert space with inner product
 $$
 \langle h_1, h_2\rangle_{\HH}:=\int_0^T\left\langle\dot h_1(s), \dot h_2(s)\right\rangle ds.
 $$

Let $\AA$ denote  the class of all $\{\FF_t\}_{0\le t \le T}$-predictable processes belonging to $\HH$ a.s..
Define for any $N \in \nn$, $$S_N:=\left\{h\in \HH; \ \int_0^T\left|\dot{h}(s)\right|^2ds\le N\right\}.$$
Consider the weak convergence topology on $\HH$, i.e.,  for any $h_n, h \in \HH, n\ge 1$, $h_n$ converges weakly to $h$ as $n \to +\infty$ if
$$
\<h_n-h,g\>_\HH \to 0, \mbox{ as } n \to +\infty,\ \forall g \in \HH.
$$
It is easy to check that $S_N$ is a compact set in $\HH$ under the weak convergence topology.
Define
$$\AA_N:=\left\{\phi\in \AA ;\ \phi(\omega)\in S_N,\ \mathbb{P}\text{-a.s.}\right\}.$$

\vskip0.3cm

We present the following result from Budhiraja et al. \cite{BDM}.

\bthm\label{thm BD}(\cite{BDM}) Let $\EE$ be a Polish space with the Borel $\sigma$-field $\BB(\EE)$. For any $\e>0$, let $\Gamma_\e$ be a measurable mapping from $C([0,T];\rr^d)$ into $\EE$.
Let $X_\e(\cdot):=\Gamma_\e(w(\cdot))$. Suppose there  exists a measurable mapping $\Gamma_0:C([0,T];\rr^d)\rightarrow \EE$ such that
\begin{itemize}
   \item[(a)]for every $N<+\infty$, the set
   $$
 \left\{\Gamma_0\left(\int_0^{\cdot}\dot h(s)ds\right);\ h\in S_N\right\}
  $$
   is a compact subset of $\EE$;
   \item[(b)] for every $N<+\infty$ and any family $\{ h^\e\}_{\e>0}\subset \AA_N$ satisfying that $h^\e$ (as $S_N$-valued random elements) converges in distribution to $h \in \AA_N$ as $\e\rightarrow 0$,
    $$\Gamma_\e\left(w(\cdot)+\frac{1}{\sqrt\e}\int_0^{\cdot}\dot h^\e(s)ds\right)\  \text{converges to} \ \Gamma_0\left(\int_0^{\cdot}\dot h(s)ds\right)$$
    in distribution as $\e\rightarrow 0$.
 \end{itemize}
Then the family $\{X_\e\}_{\e>0}$ satisfies the LDP on $\EE$ with the rate function $I$ given by
\beq\label{rate function}
I(g):=\inf_{h\in \HH;g=\Gamma_0\left(\int_0^{\cdot}\dot h(s)ds\right)}\left\{\frac12\int_0^T\left|\dot h(s)\right|^2ds\right\},\ g\in\EE,
\nneq
with the convention $\inf\emptyset=\infty$.
\nthm

\subsection{Reduction to the bounded case}
\vskip0.3cm

 Under Hypothesis \ref{hyp}, for every fixed $\e >0$, Eq.\eqref{Langevin eq time} admits a unique solution $q_\e$ in $L^k(\Omega;C([0,T];\rr^d))$. According to the proof of Theorem 3.3 in \cite{CF}, we know that the solution $q_\e$ of Eq.\eqref{Langevin eq time} can be expressed in the following form:
\beq\label{qe}
q_\e(t)= q+ \int_0^t \frac{b(q_\e(s))}{\alpha(q_\e(s))} ds + \sqrt{\e} \int_0^t \frac{\sigma(q_\e(s))}{\alpha(q_\e(s))} d w(s)+R_\e(t),
\nneq
where
\begin{align}\label{R}
R_\e(t)&:=\frac{p}{\e}\int_0^t e^{-A_\e (s)}ds -\frac{1}{\alpha(q_\e(t))} \int_0^t e^{-A_\e(t,s)} b(q_\e(s))ds\notag\\
&\quad +\int_0^t\left( \int_0^s e^{-A_\e(s,r)} b(q_\e(r)) dr \right)\frac{1}{\alpha^2(q_\e(s))} \< \nabla\alpha(q_\e(s)), \dot{q}_\e (s)\>ds\notag\\
&\quad -\frac{1}{\alpha(q_\e(t))}H_\e (t) + \int_0^t \frac{1}{\alpha^2(q_\e(s))} H_\e(s)\< \nabla\alpha(q_\e(s)), \dot{q}_\e (s)\>ds\notag\\
&=:\sum_{k=1}^5 I_\e^k (t),
\end{align}
with
$$
\begin{aligned}
&A_\e (t,s):=\frac{1}{\e^2} \int_s^t \alpha(q_\e(r))dr,\quad A_\e(t):=A_\e(t,0),\\
&H_\e(t): =\sqrt{\e} e^{-A_\e(t)} \int_0^t e^{A_\e (s)} \sigma(q_\e(s)) d w(s).
\end{aligned}
$$

We denote  the solution functional from $C([0,T];\rr^d)$ into $C([0,T];\rr^d)$ by $\mathcal G_{\e}$, i.e.,
\beq\label{G e}
\mathcal G_{\e}(w(t)):=q_\e(t),\ \forall t \in[0,T] .
\nneq
Let
\beq\label{gamma e}  X_{\e}(t):=\Gamma_\e(w(t)) :=\frac{\mathcal G_{\e}(w(t))-q_{0}(t)}{\sqrt{\e} h(\e)}, \ \forall t \in [0,T].\nneq
 Then $X_\e$ solves the following equation
\begin{align}\label{Xe}
X_\e(t)&= \frac{1}{\sqrt{\e} h(\e)} \int_0^t \left[ \frac{b(q_{0}(s)+\sqrt{\e} h(\e)X_\e(s))}{\alpha(q_{0}(s)+\sqrt{\e} h(\e)X_\e(s))}-\frac{b(q_{0}(s))}{\alpha(q_{0}(s))}\right]ds\notag\\
&\quad +\frac{1}{h(\e)} \int_0^t \frac{\sigma(q_{0}(s)+\sqrt{\e} h(\e)X_\e(s))}{\alpha(q_{0}(s)+\sqrt{\e} h(\e)X_\e(s))} dw(s)+\frac{R_\e(t)}{\sqrt{\e}h(\e)}, \quad t\in[0,T].
\end{align}

We shall prove that $\{X_\e\}_{\e >0}$ obeys an LDP on $C([0,T]; \rr^d)$ with speed function $h^2(\e)$ and the rate function $I$ given by \eqref{rate function 1}.

Since  the  family $\{q_\e\}_{\e>0}$  satisfies the  LDP in the space $C([0,T]; \rr^d)$ with the rate function $I$ and the speed function $\e^{-1}$  under Hypothesis \ref{hyp} (see Cerrai and Freidlin \cite{CF}),  there exist  some positive constants $R,C$ such that
$$
\limsup_{\e \to 0} \e \log \pp \left( \|q_\e\| \ge R \right) \le -C.
$$
Noticing \eqref{h}, we have
\beq\label{LDPbound2}
\begin{aligned}
\limsup_{\e \to 0} \frac{1}{h^2(\e)} \log \pp \left( \|q_\e\| \ge R \right)=-\infty.
\end{aligned}
\nneq

For any fixed constant $M>R$, define
$$
b^M(x):=\begin{cases}
b(x), \ &|x|< M;\\
g(x), \ &M\le|x|< M+1;\\
0, \ &|x|\ge M+1,
\end{cases}
$$
where $ g(x)$ is some infinitely differentiable function on $\rr^d$ such that $b^M(x)$ is continuous differentiable on $\rr^d$. Then for all $t\in [0,T]$, we denote
$$
\begin{aligned}
&q_0^M(t):=q+\int_0^t \frac{b^M(q_0^M(s))}{\alpha(q_0^M(s))}ds;\\
&q_\e^M(t):= q+ \int_0^t \frac{b^M(q_\e^M(s))}{\alpha(q_\e^M(s))} ds + \sqrt{\e} \int_0^t \frac{\sigma(q_\e^M(s))}{\alpha(q_\e^M(s))} d w(s)+R^M_\e(t);\\
&X_\e^M(t):= \frac{1}{\sqrt{\e} h(\e)} \int_0^t \left[ \frac{b^M(q_{0}^M(s)+\sqrt{\e} h(\e)X_\e^M(s))}{\alpha(q_{0}^M(s)+\sqrt{\e} h(\e)X_\e^M(s))}-\frac{b^M(q_{0}^M(s))}{\alpha(q_{0}^M(s))}\right]ds\\
&\quad \quad \quad \quad +\frac{1}{h(\e)} \int_0^t \frac{\sigma(q_{0}^M(s)+\sqrt{\e} h(\e)X_\e^M(s))}{\alpha(q_{0}^M(s)+\sqrt{\e} h(\e)X_\e^M(s))} dw(s)+\frac{R_\e^M(t)}{\sqrt{\e}h(\e)},
\end{aligned}
$$
where the expression of $R^M_\e(t)$ is similar to Eq.\eqref{R} with $b^M, q_\e^M$ in place of $b, q_\e$.

Notice that $\|q_0\|$ is finite by the continuity of $b$ and $\alpha$. Hence, we can choose $M$ large enough such that
$
q_0(t)=q_0^M(t),\  \mbox{for all} \ t\in[0,T].
$
Then for some $M$ large enough, by Eq.\eqref{LDPbound2}, for all $\delta>0$, we have
\begin{align}
&\limsup_{\e \to 0} \frac{1}{h^2(\e)}\log \pp(\left\|X_\e-X_\e^M\right\| >\delta)\notag\\
=&\limsup_{\e \to 0} \frac{1}{h^2(\e)}\log \pp\left(\left\|\frac{q_\e-q_\e^M}{\sqrt{\e}h(\e)}\right\| >\delta\right)\notag\\
\le &\limsup_{\e \to 0} \frac{1}{h^2(\e)}\log \pp(\left\|q_\e-q_\e^M\right\| >0)\notag\\
\le &\limsup_{\e \to 0} \frac{1}{h^2(\e)}\log \pp(\|q_\e\| \ge M)=-\infty ,
\end{align}
which means that $X_\e$ is $h^2(\e)$-exponentially equivalent to $X_\e^M$.
Hence, to prove the LDP for $\{X_\e\}_{\e>0}$ on $ C([0,T];\rr^d)$, it is enough to prove that for $\{X_\e^M\}_{\e>0}$, which is the task of the next part.

\subsection{The LDP for $\{X_\e^M\}_{\e>0}$}
In this subsection, we prove that for some fixed constant $M$ large enough , $\{X_\e^M\}_{\e >0}$ obeys an LDP on $C([0,T]; \rr^d)$ with speed function $h^2(\e)$ and the rate function $I$ given by \eqref{rate function 1}. Without loss of generality, we assume that $b$ is bounded, i.e., $|b| \le K$ for some positive constant $K$. Then $\frac{b}{\alpha}$ is also Lipschitz continuous and bounded, and by the differentiability of $\frac{b}{\alpha}$, $D(\frac{b}{\alpha})$ is also bounded. From now on, we can drop the $M$ in the notations for the sake of simplicity.

\subsubsection{Skeleton Equations}

%For any $h\in\HH$, recall the mapping  $\Gamma_0$ defined by \eqref{Gamma} that will be used to define the rate function and also used for verification of
%conditions in Theorem \ref{thm BD}.

For any $h\in \HH$, consider the deterministic equation:
\beq\label{eq skeleton}
g^h(t)= \int_0^t D\left( \frac{b(q_0(s))}{\alpha(q_0(s))}\right) g^{h}(s)ds +\int_0^t \frac{\sigma (q_0(s))}{\alpha(q_0(s))} \dot{h}(s) ds.
\nneq

\blem\label{lem skeleton} Under Hypothesis \ref{hyp}, for any $h\in \HH$, Eq.\eqref{eq skeleton} admits a
 unique  solution $g^h$ in $C([0,T];\rr^d)$, denoted by $g^h(\cdot)=:\Gamma_0\left(\int_0^\cdot \dot h(s)ds\right)$.
  Moreover, for any $N>0$, there exists some positive constant  $C(K,N,T,\alpha_0,\alpha_1)$  such that
\beq\label{eq skeleton estimate}
\sup_{h\in S_N}\left\|g^h\right\|\le C(K,N,T,\alpha_0,\alpha_1).
\nneq

\nlem

\bprf\ \ The existence and uniqueness of the solution can be proved
similarly to  the case of stochastic differential equation \eqref{eq g}, but much more simply.   \eqref{eq skeleton estimate} follows from the boundness conditions of the coefficient functions and Gronwall's inequality. Here we omit the relative proof.
\nprf

\bprop\label{Prop Gamm 0 compact}
Under Hypothesis \ref{hyp},
  for every positive number $N<+\infty$,
  the family
  $$K_N:= \left\{\Gamma_0\left(\int_0^{\cdot}\dot h(s)ds\right); h\in S_N\right\}$$
  is compact in $C([0,T];\rr^d)$.
 \nprop
\bprf\
To prove this proposition, it is sufficient to prove that the mapping $\Gamma_0$ defined in Lemma \ref{lem skeleton} is continuous from $S_N$ to $C([0,T];\rr^d)$, since  the fact that $K_N$ is compact follows from the compactness of $S_N$ under the weak topology and the continuity of the mapping $\Gamma_0$ from $S_N$ to $C([0,T];\rr^d)$.

Assume that  $h_n\to h$ weakly in $S_N$ as $n\to \infty$. We consider the following equation
\begin{align*}
 &g^{h_n}(t)-g^h(t)\\
 =& \int_0^tD\left( \frac{b(q_0(s))}{\alpha(q_0(s))}\right) \left(g^{h_n}(s)-g^h(s)\right)ds+\int_0^t \frac{\sigma (q_0(s))}{\alpha(q_0(s))}\left(\dot h_n(s)-\dot h(s)\right)ds\\
 =:&I^n_1(t)+I^n_2(t).
\end{align*}
Due to Cauchy-Schwartz inequality and the boundness of functions $\sigma,\alpha$, we know that for any $0 \le t_1 \le t_2 \le T$,
 \begin{align}\label{eq AA2}
 \left|I_2^n(t_2)-I_2^n(t_1)\right| = & \left|\int_{t_1}^{t_2}\frac{\sigma (q_0(s))}{\alpha(q_0(s))} \left(\dot h_n(s)-\dot h(s)\right)ds\right|\notag\\
 \le & \left(\int_{t_1}^{t_2} \left\|\frac{\sigma (q_0(s))}{\alpha(q_0(s))}\right\|^2_{HS}ds\right)^{\frac12}\cdot \left(\int_{t_1}^{t_2} \left|\dot h_n(s)-\dot h(s)\right|^2ds\right)^{\frac12}\notag\\
 \le &  C(K, \alpha_0)N^{\frac12}(t_2-t_1)^{\frac12}.
\end{align}
 Hence, the family of functions $\{I_2^n\}_{n\ge1}$ is equiv-continuous in $C([0,T];\mathbb R^d)$. Particularly, taking $t_1=0$, we obtain that
 \begin{align}\label{eq AA3}
\left\|I^n_2\right\|\le &  C( K, N, T,\alpha_0)<\infty,
\end{align}
where $C( K, N, T,\alpha_0)$ is independent of $n$. Thus, by the Ascoli-Arzel\'a theorem, the set $\{I_2^n\}_{n\ge1}$ is compact in $C([0,T];\rr^d)$.

 On the other hand,  for any $v\in \rr^d$, by the boundness of $\sigma/\alpha$,  we know that  the function $  \frac{\sigma(q_0)}{\alpha(q_0)}v$ belongs to  $L^2([0,T];\rr^d)$.
 Since $\dot h_n\to \dot h$ weakly in $L^2([0,T];\rr^d)$ as $n \to +\infty$, we know that
 \begin{align}\label{eq AA1}
 \left\langle I^n_2(t),v\right\rangle =\int_0^t \frac{\sigma (q_0(s))}{\alpha(q_0(s))}\left(\dot h_n(s)-\dot h(s)\right)vds\to 0, \ \ \text{as } n\rightarrow\infty.
\end{align}
Then by the compactness of $\{I_2^n\}_{n\ge 1}$, we have
\beq\label{eq I2}
\lim_{n\rightarrow\infty}\left\|I^n_2\right\|=0.
 \nneq

Set $\zeta^n(t)=\sup_{0\le s\le t}\left|g^{h_n}(s)-g^h(s)\right|$. By the boundness of $D(b/\alpha)$, we have
   \begin{align*}
   \zeta^n(t)\le C(K,\alpha_0,\alpha_1)\int_0^t \zeta^n(s) ds+\left\|I^n_2\right\|.
\end{align*}
By Gronwall's inequality and \eqref{eq I2}, we have
  $$
  \left\|g ^{h_n}-g^h\right\|\le e^{C(K,\alpha_0,\alpha_1)T}\cdot  \left\|I^n_2\right\|\to 0, \text{ as } n\to \infty,
  $$
which completes the proof.
 \nprf

\subsubsection{MDP}

For any predictable process $\dot{u}$ taking values in $L^2 ([0,T]; \rr^d)$, we denote by $q_\e^u(t)$ the solution of the following equation
\beq\label{Langevin eq u}
\begin{cases}
\e^2 \ddot{q}_\e^u (t)=b(q_\e^u(t))-\alpha(q_\e^u(t)) \dot{q}_\e^u (t) +\sqrt{\e}\sigma (q_\e^u(t)) \dot{w}(t)+\sqrt{\e} h(\e)\sigma(q_\e^u(t) )\dot{u}(t), \ t \in[0,T],\\
q_\e^u(0)=q \in \rr^d,\quad \dot{q}_\e^u(0)=\frac {p}{\e} \in \rr^d.
\end{cases}
\nneq
As is well known, for any fixed $\e >0$, $T>0$ and $k\ge 1$, this equation admits a unique solution $q_\e^u$ in $L^k(\Omega; C([0,T];\rr^d))$ as follows
$$
q_\e^{u}(t)=\GG_\e \left(w(t)+h(\e)\int_0^{t}\dot u(s)ds\right),
$$
where $\GG_\e$ is defined by \eqref{G e}.

\vskip0.3cm
%We give this proof in the appendix.
\blem\label{lem Y} Under Hypothesis \ref{hyp}, for every fixed  $N\in\mathbb{N}$ and $\e >0$, let $u^\e\in \mathcal{A}_N$ and $\Gamma_\e$ be given by \eqref{gamma e}. Then $X_\e^{u^\e}(\cdot):=\Gamma_\e\left(w(\cdot)+h(\e)\int_0^{\cdot}\dot u^\e(s)ds\right)$ is the unique solution of the following equation
\begin{align}\label{eq Y e l}
X_\e^{u^\e} (t) &= \int_0^t \frac{1}{\sqrt{\e} h(\e)} \left[ \frac{b( q_{0} (s) + \sqrt{\e} h(\e)X_\e^{u^\e}(s) )}{ \alpha(q_{0} (s) + \sqrt{\e} h(\e)X_\e^{u^\e}(s))} - \frac{b(q_{0} (s))}{\alpha(q_{0}(s))} \right] ds\notag\\
&\quad+ \int_0^t \frac{\sigma(q_{0} (s) + \sqrt{\e} h(\e)X_\e^{u^\e}(s))}{\alpha(q_{0}(s) + \sqrt{\e} h(\e)X_\e^{u^\e}(s))} \dot{u}^\e(s) ds\notag\\
&\quad+ \frac{1}{h(\e)} \int_0^t \frac{\sigma(q_{0} (s) + \sqrt{\e} h(\e)X_\e^{u^\e}(s))}{ \alpha(q_{0} (s) + \sqrt{\e} h(\e)X_\e^{u^\e}(s))} dw(s)+\frac{R_\e^{u^\e}(t)}{\sqrt{\e} h(\e)}, \quad t\in[0,T],
\end{align}
where
$$
\begin{aligned}
R_\e^{u^\e} (t)&=\frac{p}{\e}\int_0^t e^{-A_\e^{u^\e} (s)}ds -\frac{1}{\alpha(q_\e^{u^\e}(t))} \int_0^t e^{-A_\e^{u^\e}(t,s)} b(q_\e^{u^\e}(s))ds\\
&\quad +\int_0^t\left( \int_0^s e^{-A_\e^{u^\e}(s,r)} b(q_\e^{u^\e}(r)) dr \right)\frac{1}{\alpha^2(q_\e^{u^\e}(s))} \left\< \nabla\alpha(q_\e^{u^\e}(s)), \dot{q}_\e^{u^\e} (s)\right\>ds\\
&\quad -\frac{1}{\alpha(q_\e^{u^\e}(t))}H_\e^{1,u^\e} (t) + \int_0^t \frac{1}{\alpha^2(q_\e^{u^\e}(s))} H_\e^{1,u^\e}(s)\left\< \nabla\alpha(q_\e^{u^\e}(s)), \dot{q}_\e^{u^\e} (s)\right\>ds\\
&\quad -\frac{1}{\alpha(q_\e^{u^\e}(t))}H_\e^{2,u^\e} (t) + \int_0^t \frac{1}{\alpha^2(q_\e^{u^\e}(s))} H_\e^{2,u^\e}(s)\left\< \nabla\alpha(q_\e^{u^\e}(s)), \dot{q}_\e^{u^\e} (s)\right\>ds\\
&=:\sum_{k=1}^7 I_\e^{k,u^\e},
\end{aligned}
$$
with
\beq\label{notations}
\begin{aligned}
&A_\e^{u^\e} (t,s):=\frac{1}{\e^2} \int_s^t \alpha(q_\e^{u^\e}(r))dr,\quad A_\e^{u^\e}(t)=A_\e^{u^\e}(t,0),\\
&H_\e^{1,u^\e}(t):=\sqrt{\e} e^{-A_\e^{u^\e}(t)} \int_0^t e^{A_\e^{u^\e} (s)} \sigma(q_\e^{u^\e}(s)) d w(s),\\
&H_\e^{2,u^\e}(t):=\sqrt{\e}h(\e) e^{-A_\e^{u^\e}(t)} \int_0^t e^{A_\e^{u^\e} (s)} \sigma(q_\e^{u^\e}(s)) \dot{u}^\e(s)ds.
\end{aligned}
\nneq
Furthermore, there exists a positive constant $\e_0 >0$ such that for any $\e \in (0,\e_0]$,
\beq\label{L2 bounded} \ee\left[\int_0^T \left|X_\e^{u^\e} (t)\right|^2dt\right]\le C(K,N,T,\alpha_0,\alpha_1,|p|,|q|).\nneq
Moveover, we have
\beq\label{supbound}
\ee \left[\left\|X_\e^{u^\e}\right\|^2\right]  \le C(K,N,T,\alpha_0,\alpha_1,|p|,|q|).
\nneq
\nlem

\vskip0.3cm
To prove Lemma \ref{lem Y} and our main result, we present the following three lemmas.
The first lemma is similar to \cite[Lemma 3.1]{CF}.

\blem\label{LemH1}
Under Hypothesis \ref{hyp}, for any $T>0$, $k\ge 1$ and $N>0$, there exists some constant $\e_0>0$ such that for any $u^\e \in \AA_N$ and $\e \in (0,\e_0]$, we have
\beq\label{H1}
\sup_{t \in[0,T] } \ee\left[\left|H_\e^{1,u^\e}(t)\right|^k\right] \le C(k,K ,N, T,\alpha_0,\alpha_1)\left(|q|^k+|p|^k +1\right) \e^{\frac{3k}{2}}+C(k,K) \e^{\frac{k}{2}} t^{\frac{k}{2}} e^{-\frac{k\alpha_0 t}{\e^2}}.
\nneq
Moveover, we have
\beq\label{H11}
\ee \left\|H_\e^{1,u^\e}\right\|\le \sqrt{\e} C(K,N,T,\alpha_0,\alpha_1)(1+|q|+|p|).
\nneq
\nlem

\bprf Notice that Eq.\eqref{Langevin eq u} can be rewritten as the following equation: for all $t \in [0,T]$,
$$
\begin{cases}
\dot{q}_\e^{u^\e}(t)=p_\e^{u^\e}(t),\\
\e^2\dot{p}_\e^{u^\e}(t)=b(q_\e^{u^\e}(t))-\alpha(q_\e^{u^\e}(t))p_\e^{u^\e}(t)+\sqrt{\e}\sigma(q_\e^{u^\e}(t))\dot{w}(t)+
\sqrt{\e}h(\e)\sigma(q_\e^{u^\e}(t))\dot{u}^\e(t),\\
q_\e^{u^\e}(0)=q\in \rr^d, \ p_\e^{u^\e}(0)=\frac{p}{\e} \in \rr^d.
\end{cases}
$$
From the notation given in Eq.\eqref{notations}, we have
\beq\label{dotq}
\dot{q}_\e^{u^\e}(t)=p_\e^{u^\e}(t) =\frac{1}{\e}e^{-A_\e^{u^\e}(t)} p+\frac{1}{\e^2}\int_0^t e^{-A_\e^{u^\e}(t,s)}b(q_\e^{u^\e}(s))ds+ \frac{1}{\e^2}H_\e^{2,u^\e}(t)+\frac{1}{\e^2}H_\e^{1,u^\e}(t).
\nneq
Integrating with respect to $t$, we obtain that
$$
\begin{aligned}
q_\e^{u^\e}(t)&=q+\frac{1}{\e}\int_0^te^{-A_\e^{u^\e}(s)} pds+\frac{1}{\e^2}\int_0^t\int_0^s e^{-A_\e^{u^\e}(s,r)}b(q_\e^{u^\e}(r))drds\\
&\quad+\frac{1}{\e^2}\int_0^t H_\e^{2,u^\e}(s)ds+\frac{1}{\e^2}\int_0^t H_\e^{1,u^\e}(s)ds.
\end{aligned}
$$
By  Hypothesis \ref{hyp} and  Young's inequality for integral operators, we have
$$
\begin{aligned}
\left|q_\e^{u^\e}(t)\right| &\le |q|+\frac{\e}{\alpha_0} |p|+C(K,T,\alpha_0)\int_0^t\left(1+\left|q_\e^{u^\e}(s)\right|\right)ds\\
&\quad+C(K,\alpha_0)\sqrt{\e}h(\e)\int_0^t\left|\dot{u}^\e(s)\right|ds+\frac{1}{\e^2}\int_0^t \left|H_\e^{1,u^\e}(s)\right|ds\\
&\le C(K,N,T,\alpha_0)(|q|+\e |p|+\sqrt{\e}h(\e) ) +\frac{1}{\e^2}\int_0^t \left|H_\e^{1,u^\e}(s)\right|ds+C(K,T,\alpha_0)\int_0^t\left|q_\e^{u^\e}(s)\right|ds.
\end{aligned}
$$
Since $\lim_{\e \to 0}\sqrt{\e}h(\e)=0$, for $\e$ small enough, by  Gronwall's inequality,
\beq\label{qqq}
\left|q_\e^{u^\e}(t)\right|\le C(K,N,T,\alpha_0)(|q|+|p|+1) +C(K,T,\alpha_0)\frac{1}{\e^2}\int_0^t \left|H_\e^{1,u^\e}(s)\right|ds.
\nneq
Hence by the similar proof to that in \cite[Lemma 3.1]{CF}, we obtain \eqref{H1} and \eqref{H11}.

\nprf

For  $H_\e^{2,u^\e}(t)$, we have the following estimation.

\blem\label{H2}
Under Hypothesis \ref{hyp}, for any $T>0$, $k\ge 1$ and $N\in \nn$, there exists some constant $\e_0>0$ such that for any $u^\e \in \AA_N$ and $\e \in (0,\e_0]$, we have
\beq\label{H2e}
\ee \left[\left\|H_\e^{2,u^\e}\right\|^k\right]\le  C(K,N,\alpha_0)\e^{\frac{3k}{2}} h^k(\e).
\nneq
\nlem

\bprf
For any $t\in[0,T]$ and $u^\e \in \AA_N$, by the boundness of $\sigma$ and Cauchy-Schwarz inequality, we have
$$
\begin{aligned}
\left|H_\e^{2,u^\e}(t)\right| &=\left| \sqrt{\e}h(\e) e^{-A_\e^{u^\e}(t)} \int_0^t e^{A_\e^{u^\e} (s)} \sigma(q_\e^{u^\e}(s)) \dot{u}^\e(s)ds \right|\\
&\le K \sqrt{\e}h(\e)e^{-A_\e^{u^\e}(t)} \int_0^t e^{A_\e^{u^\e} (s)} \left|\dot{u}^\e(s)\right|ds\\
&\le K \sqrt{\e}h(\e)e^{-A_\e^{u^\e}(t)} \left( \int_0^t e^{2A_\e^{u^\e} (s)} ds\right)^{\frac{1}{2}} \left(\int_0^T \left|\dot{u}^\e(s)\right|^2ds \right)^{\frac{1}{2}}\\
&\le KN^{\frac{1}{2}}\sqrt{\e}h(\e)e^{-A_\e^{u^\e}(t)} \left( \int_0^t e^{2A_\e^{u^\e} (s)} ds\right)^{\frac{1}{2}}.
\end{aligned}
$$
Since $A_\e^{u^\e} (t)=\frac{1}{\e^2} \int_0^t \alpha(q_\e^{u^\e}(r))dr $, we have
$$
\begin{aligned}
\int_0^t e^{2A_\e^{u^\e} (s)} ds&=\int_0^t \frac{\e^2}{ 2\alpha(q_\e^{u^\e}(s))} de^{ \frac{2}{\e^2} \int_0^s \alpha(q_\e^{u^\e} (r))dr}\\
&\le\frac{\e^2}{ 2 \alpha_0} \int_0^t de^{ \frac{2}{\e^2} \int_0^s \alpha(q_\e^{u^\e} (r))dr}\\
&=\frac{\e^2}{ 2 \alpha_0}\left(e^{2A_\e^{u^\e}(t)}-1\right).
\end{aligned}
$$
Hence
$$
\begin{aligned}
\left|H_\e^{2,u^\e}(t)\right| &\le KN^{\frac{1}{2}}\frac{\e^\frac{3}{2} h(\e)}{ \sqrt{2 \alpha_0}}e^{-A_\e^{u^\e}(t)}\left(e^{2A_\e^{u^\e}(t)}-1\right)^\frac{1}{2}\\
&\le C(K,N,\alpha_0) \e^{\frac{3}{2}} h(\e) e^{-A_\e^{u^\e}(t)}e^{A_\e^{u^\e}(t)}\\
&=C(K,N,\alpha_0)\e^{\frac{3}{2}} h(\e),
\end{aligned}
$$
and furthermore
$$
\ee\left[\left\|H_\e^{2,u^\e}\right\|^k\right]\le  C(K,N,\alpha_0)\e^{\frac{3k}{2}} h^k(\e),
$$
which completes the proof.
\nprf

\blem\label{LemR}
Under Hypothesis \ref{hyp}, for any $T>0$ and any $u^\e \in \AA_N$, we have
\beq\label{R1}
\ee\left\|\frac{R_\e}{\sqrt{\e}h(\e)}\right\|\to 0, \ \ \ \text{  as } \ \  \e \rightarrow 0.
\nneq
Moreover, we have
\beq\label{R2}
\ee\left[ \left\|\frac{R_\e}{\sqrt{\e}h(\e)}\right\|^2\right]\to 0, \ \ \ \text{  as } \ \  \e \rightarrow 0.
\nneq
\nlem

\bprf
Similarly to the proof \cite[(3.17)]{CF},  under Hypothesis \ref{hyp},  we have
\beq \label{Y41}
\ee \left\|\frac{\sum_{k=1}^5 I_\e^{k,u^\e}}{\sqrt{\e}h(\e)}\right\| \le \frac{1}{h(\e)}C(K,N,T,\alpha_0,\alpha_1,|p|,|q|) \to 0, \  \mbox{ as } \e \rightarrow 0.
\nneq
Next, we will estimate   $\ee \left\|\frac{I_\e^{6,u^\e}}{\sqrt{\e} h(\e)}\right\|$ and $\ee \left\|\frac{I_\e^{7,u^\e}}{\sqrt{\e} h(\e)}\right\|$.
By Lemma \ref{H2}, we have
\begin{align}\label{I6}
\ee \left\|\frac{I_\e^{6,u^\e}}{\sqrt{\e} h(\e)}\right\|
 \le  \frac{1}{\sqrt{\e} h(\e)\alpha_0} \ee \left\|H_\e^{2,u^\e}\right\|
 \le  \e  C(K,N,\alpha_0) \to 0,  \ \mbox { as } \e \rightarrow 0.
\end{align}
By Cauchy-Schwarz inequality, we have
$$
\begin{aligned}
\ee \left\|\frac{I_\e^{7,u^\e}}{\sqrt{\e} h(\e)}\right\| &\le \frac{C(K,\alpha_0)}{ \sqrt{\e} h(\e)} \ee \left[\sup_{t\in[0,T]} \int_0^t \left|H_\e^{2,u^\e}(s)\right|\cdot \left|\dot{q}_\e^{u^\e}(s)\right|ds  \right]\\
&\le \frac{C(K,\alpha_0)}{ \sqrt{\e} h(\e)} \left[ \int_0^T \ee\left[\left|H_\e^{2,u^\e}(s)\right|^2 \right]ds \right]^{\frac{1}{2}}\cdot\left[\int_0^T \ee\left[\left|\dot{q}_\e^{u^\e}(s)\right|^2\right]ds  \right]^{\frac{1}{2}}.
\end{aligned}
$$
By \eqref{qqq}, we have for all $\e>0$ small enough,
$$
\int_0^T \left|\dot{q}_\e^{u^\e}(s)\right|^2 ds \le C(K,N,T,\alpha_0,|p|,|q|)+\frac{C(K,T,\alpha_0)}{\e^4}\int_0^T \left|H_\e^{1,u^\e}(s)\right|^2 ds.
$$
Hence, by \eqref{H1} and Lemma \ref{H2}, we have
\begin{align}\label{I7}
&\quad\ee\left\|\frac{I_\e^{7,u^\e}}{\sqrt{\e} h(\e)}\right\|\notag\\
&\le \frac{C(K,N,T,\alpha_0,|p|,|q|)}{\sqrt{\e} h(\e)}\left[ \left( \int_0^T \ee\left[\left|H_\e^{2,u^\e}(s)\right|^2\right] ds\right)^{\frac{1}{2}}\right]\notag\\
&\quad+ \frac{C(K,N,T,\alpha_0)}{\e^{\frac{5}{2}} h(\e)}\left( \int_0^T \ee\left[\left|H_\e^{2,u^\e}(s)\right|^2\right] ds\right)^{\frac{1}{2}}\cdot\left( \int_0^T \ee\left[\left|H_\e^{1,u^\e}(s)\right|^2\right] ds\right)^{\frac{1}{2}} \notag\\
&\le \sqrt{\e} C(K,N,T,\alpha_0,\alpha_1,|p|,|q|) \to 0,\ \ \mbox{ as } \e \rightarrow 0.
\end{align}

This together with  \eqref{Y41} and \eqref{I6} implies \eqref{R1}.

\eqref{R2} can be easily obtained by applying the similar estimation process for
$$ \ee\left[\left\|\frac{I_\e^{i,u^\e}}{\sqrt{\e} h(\e)}\right\|^2\right],\ i=1,2,3,\cdots,7, $$
as given above. Hence we omit the proof.
\nprf

Now we prove Lemma \ref{lem Y}.
\bprf[{\bf The proof of Lemma \ref{lem Y}}]
For any $\e >0$ and $u^\e \in \AA_N$, define
$$
d\qq^{u^\e} := \exp\left\{ -h(\e) \int_0^t \dot{u}^\e(s)dw(s)-\frac{ h^2(\e)}{2} \int_0^t \left|\dot{u}^{\e}(s)\right|^2 ds\right\}d\pp.
$$
Since $\frac{d\qq^{u^\e}}{d\pp}$ is an exponential martingale, $\qq^{u^\e}$ is a probability measure on $\Omega$. By Girsanov theorem, the process
$$
\tilde{w}^{\e}(t)=w(t)+h(\e)\int_0^t\dot{u}^\e(s)ds
$$
is a $\rr^d$-valued Wiener process under the probability measure $\qq^{u^\e}$. Rewriting Eq.\eqref{eq Y e l} with $\tilde{w}^{\e}(t)$, we obtain Eq.\eqref{Xe} with $\tilde{w}^{\e}(t)$ in place of $w(t)$. Let $X_\e^{u^\e}$ be the unique solution of Eq.\eqref{Xe} with $\tilde{w}^{\e}(t)$ on the space
$(\Omega,\FF,\qq^{u^\e})$. Then $X_\e^{u^\e}$ satisfies Eq.\eqref{eq Y e l}, $\qq^{u^\e}$-a.s.. By the equivalence of probability measures, $X_\e^{u^\e}$ satisfies Eq.\eqref{eq Y e l}, $\pp$-a.s..

Now we prove \eqref{L2 bounded}. By \eqref{R2}, there exists some constant $\e_0>0$ such that for any $\e \in(0,\e_0]$,
\beq\label{R3}
\ee \left[  \left\|\frac{R_\e^{u^\e}}{\sqrt{\e}h(\e)}\right\|^2 \right] \le C(K,N,T,\alpha_0,\alpha_1,|p|,|q|).
\nneq
Notice that $b/\alpha$ is Lipschitz continuous and $\sigma/\alpha$ is bounded, then we have
\begin{align}\label{Xbound}
\left|X_\e^{u^\e} (t)\right|^2 &\le C(K,\alpha_0,\alpha_1)\int_0^t \left|X_\e^{u^\e}(s)\right|^2 ds+C(K,N,T,\alpha_0)\notag\\
&\quad+ \frac{C(K,\alpha_0)}{h^2(\e)}w^2(t)+C\left|\frac{R_\e^{u^\e}(t)}{\sqrt{\e} h(\e)}\right|^2.
\end{align}
Hence by \eqref{h} and \eqref{R3}, for any $\e \in(0,\e_0]$, taking expectation in both sides in \eqref{Xbound}, we have
$$
\ee \left[\left|X_\e^{u^\e} (t)\right|^2\right] \le C(K,\alpha_0,\alpha_1)\int_0^T \ee  \left[\left|X_\e^{u^\e}(s)\right|^2\right] ds+C(K,N,T,\alpha_0,\alpha_1,|p|,|q|).
$$
By Gronwall's inequality, we get
\beq\label{supX}
\ee \left[\left|X_\e^{u^\e} (t)\right|^2\right] \le C(K,N,T,\alpha_0,\alpha_1,|p|,|q|),
\nneq
then by Fubini's theorem,
\beq\label{intX}
\ee \left[\int_0^T \left|X_\e^{u^\e} (s)\right|^2ds\right] \le C(K,N,T,\alpha_0,\alpha_1,|p|,|q|).
\nneq
First taking supremum with respect to $t\in [0,T]$ in \eqref{Xbound}, and then taking expectation in both sides, for any $\e \in(0,\e_0]$, by BDG inequality, \eqref{h}, \eqref{R3} and \eqref{intX}, we obtain that
$$
\begin{aligned}
\ee\left[\left\|X_\e^{u^\e} \right\|^2\right] &\le C(K,\alpha_0,\alpha_1)\ee\left[\int_0^T \left|X_\e^{u^\e}(s)\right|^2 ds\right]+C(K,N,T,\alpha_0,\alpha_1,|p|,|q|)\\
&\le C(K,N,T,\alpha_0,\alpha_1,|p|,|q|),
\end{aligned}
$$
which completes the proof.
\nprf

\bprop\label{P1}
Under Hypothesis \ref{hyp}, for every fixed $N \in \nn$, let $\{u^\e\}_{\e >0}$ be a family of processes in  $\AA_N $ that converges in distribution to some $u \in \AA_N $ as $\e \to 0$, as random variables taking values in the space $S_N$, endowed with the weak topology. Then
$$\Gamma_{\e} \left( w(\cdot)+h(\e )\int_0^\cdot \dot{u}^\e(s) ds\right) \to \Gamma_0\left( \int_0^\cdot \dot{u}(s) ds \right),$$
in distribution in $C([0,T]; \rr^d)$ as $\e \to 0$.
 \nprop

\bprf
By the Skorokhod representation theorem, there exists a probability basis $(\bar\Omega,\bar\FF,(\bar\FF_t),\bar\pp)$, and  on this basis, a Brownian motion $\bar w$ and a family of $\bar \FF_t$-predictable processes $\{\bar u^\e\}_{\e>0}, \bar u$   taking values in $S_N$, $\bar\pp$-a.s., such that the joint law of $(u^\e,u, w)$ under $\pp$ coincides with that of  $(\bar u^\e, \bar u, \bar w)$ under $\bar\pp$ and
$$
\lim_{\e\rightarrow0}\langle\bar u^\e-\bar u, g \rangle_\HH=0, \ \ \forall g\in\HH, \ \bar \pp\mbox{-}a.s..  \ \
$$

Let $\bar X_{\e}^{ \bar{u}^{\e}}$ be the solution of a similar equation to \eqref{eq Y e l} with $u^\e$ replaced by $\bar u^\e$ and $w$ by $\bar w$,  and let $\bar X^{\bar u}$ be the solution of a similar equation to \eqref{eq skeleton} with $h$ replaced by $ \bar u$.
Thus, to prove this proposition, it is sufficient to prove that
\beq\label{eq b conv}
\lim_{\e\rightarrow0}\left \|\bar X_{\e}^{ \bar{u}^{\e}}-\bar X^{\bar u} \right\|=0, \quad \text{  in probability}.
\nneq

From now on, we drop the bars in the notation for the sake of simplicity.

Notice that, for any $ t \in[0,T] $,
\begin{align}\label{XeX}
&\quad X_\e^{u^\e}(t) - X^u(t)\notag\\
&=\int_0^t \left\{\frac{1}{\sqrt{\e} h(\e)} \left[ \frac{b( q_{0} (s) + \sqrt{\e} h(\e)X_\e^{u^\e}(s) )}{ \alpha(q_{0} (s) + \sqrt{\e} h(\e)X_\e^{u^\e}(s))} - \frac{b(q_{0} (s))}{\alpha (q_{0}(s))} \right]-D \left( \frac{b(q_0(s))}{\alpha(q_0(s))} \right) X^u(s)  \right\} ds\notag\\
&\quad+ \int_0^t \left[\frac{\sigma(q_{0} (s) + \sqrt{\e} h(\e)X_\e^{u^\e}(s))}{\alpha(q_{0} (s) + \sqrt{\e} h(\e)X_\e^{u^\e}(s))} \dot{u}^\e(s)-\frac{\sigma(q_0(s))}{\alpha(q_0(s))} \dot{u}(s) \right] ds\notag\\
&\quad+ \frac{1}{h(\e)} \int_0^t \frac{\sigma(q_{0} (s) + \sqrt{\e} h(\e)X_\e^{u^\e}(s))}{ \alpha(q_{0} (s) + \sqrt{\e} h(\e)X_\e^{u^\e}(s))} dw(s) +\frac{R_\e^{u^\e}(t)}{\sqrt{\e} h(\e)}\notag\\
&=: \sum_{k=1}^4 Y_\e^{k,u^\e} (t).
\end{align}
We shall prove this proposition in the following four steps.

{\bf Step 1:} For the first term $Y_\e^{1,u^\e}$,   denote $ x_\e(t):=\sqrt{\e} h(\e)X_\e^{u^\e}(t)$, by Taylor's formula, there exists a random variable $\eta_\e$ taking values in $(0,1)$ such that
$$
\begin{aligned}
&\quad \left|Y_\e^{1,u^\e} (t)\right|\\
&= \left|\int_0^t  \left[D\left( \frac{b( q_{0} (s) + \eta_\e(s)x_\e(s) )}{ \alpha(q_{0} (s) + \eta_\e(s)x_\e(s))}\right) X_\e^{u^\e}(s)-D \left( \frac{b(q_0(s))}{\alpha(q_0(s))} \right) X^u(s)  \right] ds\right|\\
 &\le \left|\int_0^t D \left( \frac{b( q_{0} (s) + \eta_\e(s)x_\e(s) )}{ \alpha(q_{0} (s) + \eta_\e(s)x_\e(s))}\right)\cdot \left( X_\e^{u^\e}(s)-X^u(s)\right) ds\right|\\
&\quad +\left|\int_0^t \left[D \left( \frac{b( q_{0} (s) + \eta_\e(s)x_\e(s) )}{ \alpha(q_{0} (s) + \eta_\e(s)x_\e(s))}\right) - D \left( \frac{b(q_0(s))}{\alpha(q_0(s))} \right)\right]\cdot X^u(s) ds\right |\\
&=:y_\e^{11}(t)+y_\e^{12}(t).
\end{aligned}
$$

For the first term $ y_\e^{11}$, by the boundness of $D\left(\frac{b}{\alpha}\right)$,    we have
\beq\label{Y11}
y_\e^{11}(t) \le C(K,\alpha_0,\alpha_1) \int_0^t \left|X_\e^{u^\e}(s)-X^u(s) \right|ds.
\nneq

Next we deal with the second term $y_\e^{12}$.
For each $R>\|q_0\|$ and $\rho\in(0,1)$,
set
$$
\eta_{R,\rho}:=\sup_{|x|\le R, |y|\le R,  |x-y|\le \rho} \left|D\left(\frac{b}{\alpha}\right)(x)-D\left(\frac{b}{\alpha}\right)(y) \right|.
$$
Then by the continuous differentiability of $\frac{b}{\alpha}$, we know that for any fixed $R>0$,
$$
\lim_{\rho\rightarrow0}\eta_{R,\rho}=0.
$$
Since $\sqrt{\e} h(\e) \to 0$ as $\e \to 0$, there exists some $\e_0>0$ small enough such that for all $0<\e\le \e_0 $,
$$
\sup_{\|q_0\|\le R, \sqrt{\e}h(\e) \|X_{\e}^{u^\e}\|\le \rho }\left\|\left(D\left(\frac{b}{\alpha}\right)(q_0+\eta_{\e}\sqrt{\e}h(\e)X_{\e}^{u^\e} )-D\left(\frac{b}{\alpha}\right)(q_0) \right)X^u\right\|\le \eta _{R+1,\rho}\|X^u\|
$$
for any $\rho \in (0,1)$.

Thus, we obtain that for any $r>0, R>\|q_0\|$,
\begin{align}\label{pp}
&\pp\left(\left\|y_{\e}^{12}\right\|>r \right)\notag\\
\le & \pp\left(\sqrt{\e}h(\e)\left\|X_{\e}^{u^\e} \right\|> \rho \right)+  \pp\left( \eta _{R+1,\rho}\left\|X^u\right\|>\frac{r}{T} \right)\notag\\
\le & \frac{\e h^2(\e)}{\rho^2}\ee\left[\left\|X_{\e}^{u^\e}\right\|^2\right]+\frac{\eta _{R+1,\rho}^2T^2}{r^2}\ee\left[\left\|X^u\right\|^2\right].
\end{align}
By \eqref{eq skeleton estimate} and \eqref{supbound}, letting $\e\rightarrow 0$ and then $\rho\rightarrow0$ in \eqref{pp}, we can prove that
\beq\label{Y120}
\lim_{\e\rightarrow0}\pp\left(\left\|y_{\e}^{12}\right\|>r \right)=0, \ \ \ \ \text{for any } r>0.
\nneq

{\bf Step 2:}
For the second term $Y_\e^{2,u^\e}$ we have
$$
\begin{aligned}
&  \left|Y_\e^{2,u^\e}(t)\right|\\
 \le &\left|\int_0^t  \frac{\sigma (q_{0} (s) + \sqrt{\e} h(\e)X_\e^{u^\e}(s))}{\alpha (q_{0} (s) + \sqrt{\e} h(\e)X_\e^{u^\e}(s))} \left( \dot{u}^\e(s)-\dot{u}(s)\right) ds\right|\\
& +\left|\int_0^t \left[\frac{\sigma(q_{0} (s) + \sqrt{\e} h(\e)X_\e^{u^\e}(s))}{\alpha(q_{0} (s) + \sqrt{\e} h(\e)X_\e^{u^\e}(s))} -\frac{\sigma(q_0(s))}{\alpha(q_0(s))}\right]  \dot{u}(s)  ds\right|\\
=:& \left|Y_\e^{2,u^\e,1}(t)\right|+\left|Y_\e^{2,u^\e,2}(t)\right|.
\\
\end{aligned}
$$

Using the same argument as that in the proof of \eqref{eq I2},  we obtain that
\beq\label{Ylim}
\lim_{\e\rightarrow 0}\left\|Y_\e^{2,u^\e,1}\right\|=0, \ \text{ a.s.}.
\nneq
Since $ \left\|Y_\e^{2,u^\e,1}\right\| \le C(K,N,T,\alpha_0)$, by the dominated convergence theorem, Eq.\eqref{Ylim} implies  that
$$
\lim_{\e\rightarrow 0}\ee\left\|Y_\e^{2,u^\e,1}\right\|= 0.
$$

Due to the Lipschitz continuity of $\sigma/\alpha$,  we have
\beq\label{Y21}
\begin{aligned}
\left\|Y_\e^{2,u^\e,2}\right\|\le  C(K,\alpha_0,\alpha_1)\int_0^T   \sqrt{\e}h(\e)\left|X_\e^{u^\e}(t)\right|\cdot\left|\dot{u}(t)\right|ds.
\end{aligned}
\nneq
By \eqref{L2 bounded} and H\"{o}lder's inequality, we get
$$
\ee\left[\int_0^T \left|X_\e^{u^\e}(t)\right|\cdot\left|\dot{u}(t)\right|dt\right] \le C(K,N,T,\alpha_0,\alpha_1,|p|,|q|).
$$
Hence by \eqref{h}, we obtain that
\beq\label{Y22}
\ee\left\|Y_\e^{2,u^\e}\right\|\to 0,   \ \ \text{as } \e \rightarrow 0.
\nneq

{\bf Step 3:}
For the third term $Y_\e^{3,u^\e}$,  by BDG inequality and \eqref{h},  we have
\begin{align}\label{Y31}
\ee \left\|Y_\e^{3,u^\e}\right\| &=\frac{1}{h(\e)}\ee\left[\sup_{t\in[0,T]} \left| \int_0^t \frac{\sigma(q_{0} (s) + \sqrt{\e} h(\e)X_\e^{u^\e}(s))}{ \alpha(q_{0} (s) + \sqrt{\e} h(\e)X_\e^{u^\e}(s))} dw(s)\right|\right]\notag\\
& \le \frac{C}{h(\e)} \ee \left(\int_0^T\left\| \frac{(\sigma*\sigma^T)(q_{0} (s) + \sqrt{\e} h(\e)X_\e^{u^\e}(s))}{ \alpha^2(q_{0} (s) + \sqrt{\e} h(\e)X_\e^{u^\e}(s))}\right\|_{HS} ds\right)^{\frac{1}{2}}\notag\\
&\le\frac{C(K,T,\alpha_0) }{h(\e)} \to 0,\  \mbox{ as } \e \rightarrow 0.
\end{align}

{\bf Step 4:}
For the last term $Y_\e^{4,u^\e}$, by Lemma \ref{LemR}, we have
\beq\label{Y42}
\ee\left\|Y_\e^{4,u^\e}\right\|\to 0, \ \ \ \text{  as } \ \  \e \rightarrow 0.
%Y_\e^{4,u^\e} (t)
\nneq
By Eq.\eqref{XeX} and \eqref{Y11}, we obtain that
\begin{align}
&\quad\sup_{0\le s\le t}\left|X_\e^{u^\e}(s)-X^u(s)\right| \notag\\
&\le C(K,\alpha_0,\alpha_1) \int_0^t \sup_{0\le v\le s}\left|X_\e^{u^\e}(v)-X^u(v) \right|ds+
\sup_{0\le s\le t}y_\e^{12}(s)\notag \\
&\quad+\sup_{{0\le s\le t}}\left|Y_\e^{2,u^\e}(s)\right|+\sup_{{0\le s\le t}}\left|Y_\e^{3,u^\e}(s)\right|+\sup_{{0\le s\le t}}\left|Y_\e^{4,u^\e}(s)\right|.
\end{align}
Using Gronwall's inequality, we have that
$$
\left\|X_\e^{u^\e}-X^u\right\|\le C\left(\left\|y_\e^{12}\right\|+\sum_{l=2,3,4}
\left\|Y_\e^{l,u^\e}\right\|\right).
$$
This, together with  \eqref{Y120}, \eqref{Y22}, \eqref{Y31} and  \eqref{Y42},   implies that
$$
\lim_{\e\rightarrow0} \left\|X_\e^{u^\e}-X^u\right\|=0, \quad \text{in probability},
$$
which completes the proof.
\nprf

According to Theorem \ref{thm BD}, the MDP of $\{X_\e^M\}_{\e >0}$ follows from Proposition \ref{Prop Gamm 0 compact} and Proposition \ref{P1}, which completes the proof of our main result Theorem \ref{MDP}.

\vskip0.3cm

\noindent{\bf Acknowledgements:} We thank the anonymous referees for their valuable comments and suggestions which help us improve the quality of this paper. Liu W. is supported by Natural Science Foundation of China (11571262, 11731009).

\end{document}